\pdfoutput=1
\documentclass[11pt]{article}
\usepackage[left=1in,right=1in,top=1in,bottom=1in]{geometry}
\usepackage{times}
\usepackage{expl3}
\usepackage{cite}
\usepackage[table]{xcolor}
\usepackage{multirow}
\usepackage{stackengine} 
\usepackage{hhline}
\usepackage{lipsum}
\usepackage{titlesec}
\usepackage{wrapfig}
\usepackage{enumerate}
\usepackage{epsfig}
\usepackage{amsmath}
\usepackage{tabularx}
\usepackage{array}
\usepackage{booktabs}
\usepackage{enumitem}
\usepackage{bbm}
\usepackage{calc}
\usepackage{graphicx}
\usepackage{amsmath}
\usepackage[title]{appendix}
\usepackage{amssymb}
\usepackage{epstopdf}
\usepackage{boldline}
\usepackage{arydshln}
\usepackage{calligra}
\usepackage{bm}
\usepackage{url}
\usepackage{blindtext}
\usepackage{accents}

\newcommand{\define}{\stackrel{\mbox{\tiny def}}{=}}

\newtheorem{theorem}{Theorem}
\newtheorem{corollary}{Corollary}

\newtheorem{example}{Example}
\newtheorem{remark}{Remark}
\usepackage{mathtools}
\usepackage{epstopdf}
\usepackage{balance}
\usepackage{thmtools}
\usepackage{thm-restate}
\usepackage{hyperref}
\usepackage{cleveref}
\usepackage[mathscr]{euscript}

\usepackage[ruled,vlined]{algorithm2e}
\include{pythonlisting}
\newcommand{\ostar}{\mathbin{\mathpalette\make@circled\star}}

\makeatletter
\newcommand{\removelatexerror}{\let\@latex@error\@gobble}
\makeatother
\makeatletter
\newcommand*{\rom}[1]{\expandafter\@slowromancap\romannumeral #1@}
\makeatother

\ExplSyntaxOn
\newcommand\latinabbrev[1]{
  \peek_meaning:NTF . {
    #1\@}%
  { \peek_catcode:NTF a {
      #1.\@ }%
    {#1.\@}}}
\ExplSyntaxOff

\titleclass{\subsubsubsection}{straight}[\subsubsection]

\begin{document}
\vspace{1cm}
\title{Tensor Unfolding Characterization}\vspace{1.8cm}
\author{Shih~Yu~Chang
\thanks{Shih Yu Chang is with the Department of Applied Data Science,
San Jose State University, San Jose, CA, U. S. A. (e-mail: {\tt
shihyu.chang@sjsu.edu}). 
           }}

\maketitle

\begin{abstract}
Tensors play a pivotal role in the realms of science and engineering, particularly in the realms of data analysis, machine learning, and computational mathematics. The process of unfolding a tensor into matrices, commonly known as tensor unfolding or matricization, serves as a valuable technique for simplifying the representation of tensors with higher orders. In this study, we initially derive unfolded matrices from a specified tensor over a B{'e}zout ring using a matrix equivalence relation. We proceed to elucidate the relationships between eigenvalues and eigenvectors within these unfolded matrices. Additionally, we employ the localization approach outlined by Gerstein to ascertain the count of distinct matrix equivalence classes present among the unfolded matrices.
\end{abstract}


\section{Introduction}\label{sec: Introduction}

Tensors play a crucial role in the fields of science and engineering, especially in the context of data analysis, machine learning, and computational mathematics. We discuss how the following fields are affected by tensors. Tensors are employed in simulations and modeling in physics and engineering. They help represent physical quantities that vary in multiple dimensions. For example, stress and strain tensors are used in materials science and structural engineering simulations~\cite{harrigan1984characterization}. Tensors are fundamental to the field of machine learning, particularly in deep learning. Neural networks, the backbone of deep learning models, use tensors to store and process data. Tensors facilitate the flow of information through the network, enabling the training and inference processes~\cite{chang2023personalized,chang2023TEKF}. In engineering, signals are often represented as tensors. For instance, in the field of image processing, tensors are used to represent and process signals, allowing for tasks such as image denoising, compression, and feature extraction~\cite{chang2023TLS}. Tensors are used in numerical methods and simulations to solve partial differential equations and other mathematical problems. They enable efficient representation and manipulation of large datasets in numerical computations~\cite{khoromskij2018tensor}. In summary, tensors serve as a foundational concept in various scientific and engineering disciplines, providing a powerful mathematical framework for representing and analyzing complex data structures. Their impact is particularly pronounced in fields where data is multi-dimensional and intricate computations are required.

Unfolding a tensor into matrices, often referred to as tensor unfolding or matricization, is a technique used to simplify the representation of higher-order tensors. This process is particularly common in the context of tensor decompositions, where the goal is to express a high-dimensional tensor as a sum of simpler components. Unfolding helps in applying well-established linear algebraic methods that operate on matrices, making it easier to analyze and manipulate the data~\cite{liang2019further}. An unfolded tensor refers to a matrix created by systematically rearranging the entries of the tensor into a two-dimensional array. Let $\mathfrak{R}$ be a B{\'e}zout ring. Consider a tensor space in $\mathfrak{R}^{I_1 \times \cdots \times I_M \times I_1 \times \cdots \times I_M}$ and a matrix space in $\mathfrak{R}^{\mathbb{I}_M \times \mathbb{I}_M}$, we define an unfolding map, represented by two tuples $(\psi_r, \psi_c)$, as follows:
\begin{eqnarray}\label{eq:unfolding map def}
(\psi_r, \psi_c): \mathfrak{R}^{I_1 \times \cdots \times I_M \times I_1 \times \cdots \times I_M}&\rightarrow&\mathfrak{R}^{\mathbb{I}_M \times \mathbb{I}_M}\nonumber \\
\mathcal{B}=[b_{\bm{i}_M,\bm{j}_M}]&\rightarrow& \bm{B}=[b_{\psi_r(\bm{i}_M,\bm{I}_M),\psi_c(\bm{j}_M,\bm{I}_M)}],
\end{eqnarray}
where $\psi_r(\bm{i}_M,\bm{I}_M) \in [\mathbb{I}_M]$ and $\psi_c(\bm{j}_M,\bm{I}_M) \in [\mathbb{I}_M]$ are two index mapping functions for the row part and the column part, respectively. For example, we can have the following mapping relations for $\psi_r(\bm{i}_M,\bm{I}_M)$ and $\psi_c(\bm{j}_M,\bm{I}_M)$ expressed by
\begin{eqnarray}\label{eq:unfolding map exp r}
\psi_r(\bm{i}_M,\bm{I}_M)&\define&i_1 + \sum\limits_{k=2}^M (i_k-1)\prod\limits_{\ell=1}^{k-1} I_\ell,
\end{eqnarray}
and
\begin{eqnarray}\label{eq:unfolding map exp c}
\psi_c(\bm{j}_M,\bm{I}_M)&\define&j_1 + \sum\limits_{k=2}^M (j_k-1)\prod\limits_{\ell=1}^{k-1} I_\ell.
\end{eqnarray}
Many works, for example~\cite{liang2019further}, use the map $(\psi_r,\psi_c)$ provided by Eqs~\eqref{eq:unfolding map exp r} and~\eqref{eq:unfolding map exp c} to relate tensor spectrum theorem and tensor SVD for the tensor $\mathcal{B}$ by the matrix spectrum theorem and matrix SVD for the matrix $\bm{B}$. In this work, we will consider two situations for $(\psi_r, \psi_c)$. If $\psi_r=\psi_c$, the eigenvalues and eigenvectors of the tensor $\mathcal{B}$ will be characterized by the eigenvalues and eigenvectors of the matrix $\bm{B}$ with respect to different mapping $\psi_r$. This is the main topic investigated in Section~\ref{sec: Row/Column Same Mapping}. On the other hand, if $\psi_r\neq\psi_c$, the equivalent classes with respect to different $\psi_r$ and $\psi_c$ will be studied in Section~\ref{sec: Row/Column Different Mapping}.

In this work, we first fomulate unfolded matrices from a given tensor over a B{\'e}zout ring via matrix equivalence relation, and characterize the eigenvalues and eigenvectors relationships among unfolded matrices. We also adopt the localization approach from~\cite{gerstein1977local} to determine the number of different matrix equivalence classes among unfolded matrices.

In Section~\ref{sec: Row/Column Same Mapping}, the spectrum relation among various unfolded matrices under $\psi_r=\psi_c$  is discussed. In Section~\ref{sec: Row/Column Different Mapping}, we consider the number of equivalence classes and spectrum for more general unfolded matrices with $\psi_r\neq\psi_c$. 

\textbf{Nomenclature:} $\mathcal{A}$, $\bm{A}$ and $\bm{a}$ represent tensor, matrix and vector, respectively. $\mathfrak{R}$ used to represent a B{\'e}zout ring. $\mathbb{R}$, $\mathbb{Z}$ and $\mathbb{N}$ represent real numbers, integers and natural numbers, respectively. $\mathbb{I}_M$ is defined as $\prod\limits_{k=1}^M I_k$, where $I_k \in \mathbb{N}$. $\mathbb{F}$ represents a field. $[N]$ denotes the set $\{1,2,\cdots, N\}$. $\bm{i}_M$ represents $M$ tuple indices $[i_1,\cdots i_M]$ and $\bm{I}_M$ represents $M$ tuple dimensions $[I_1,\cdots,I_M]$. 

\section{Row/Column Same Mapping}\label{sec: Row/Column Same Mapping}

Given a tensor $\mathcal{B} \in \mathfrak{R}^{I_1 \times \cdots \times I_M \times I_1 \times \cdots \times I_M}$, under the unfolding map $(\psi_r, \psi_c)$, we have
\begin{eqnarray}\label{eq1:sec: Row/Column Same Mapping}
\mathcal{B}=[b_{\bm{i}_M,\bm{j}_M}]&\rightarrow& \bm{B}=[b_{\psi_r(\bm{i}_M,\bm{I}_M),\psi_c(\bm{j}_M,\bm{I}_M)}].
\end{eqnarray}
If we consider another unfolding map $(\psi'_r, \psi'_c)$ for the same tensor $\mathcal{B} \in \mathfrak{R}^{I_1 \times \cdots \times I_M \times I_1 \times \cdots \times I_M}$, we have
\begin{eqnarray}\label{eq2:sec: Row/Column Same Mapping}
\mathcal{B}=[b_{\bm{i}_M,\bm{j}_M}]&\rightarrow& \bm{B}'=[b_{\psi'_r(\bm{i}_M,\bm{I}_M),\psi'_c(\bm{j}_M,\bm{I}_M)}].
\end{eqnarray}

Because both maps $(\psi_r, \psi_c)$ and $(\psi'_r, \psi'_c)$ are unfolding maps with $\psi_r=\psi_c$ and $\psi'_r=\psi'_c$, there exists a permutation matrix $\bm{P} \in \mathfrak{R}^{\mathbb{I}_M \times \mathbb{I}_M}$ such that 
\begin{eqnarray}\label{eq3:sec: Row/Column Same Mapping}
\bm{B}'=\bm{P}^{-1}\bm{B}\bm{P},
\end{eqnarray}
where $\bm{P}$ is the permutation matrix from the index map given by $\psi_r$ to the index map given by $\psi'_r$. We have the following theorem about various unfolding maps for a given tensor $\mathcal{B} \in \mathfrak{R}^{I_1 \times \cdots \times I_M \times I_1 \times \cdots \times I_M}$.
\begin{theorem}\label{thm: unfolding maps spectrum and eigenvectors}
Given a tensor $\mathcal{B}=[b_{\bm{i}_M,\bm{j}_M}] \in \mathfrak{R}^{I_1 \times \cdots \times I_M \times I_1 \times \cdots \times I_M}$ and two unfolding maps $(\psi_r, \psi_c)$ and $(\psi'_r, \psi'_c)$ with $\psi_r=\psi_c$ and $\psi'_r=\psi'_c$, we have 
\begin{eqnarray}
\bm{B}=[b_{\psi_r(\bm{i}_M,\bm{I}_M),\psi_c(\bm{j}_M,\bm{I}_M)}],
\end{eqnarray}
and
\begin{eqnarray}
\bm{B}'=[b_{\psi'_r(\bm{i}_M,\bm{I}_M),\psi'_c(\bm{j}_M,\bm{I}_M)}].
\end{eqnarray}
Then, $\bm{B}'$ and $\bm{B}$ have the same eigenvalues, and the eigenvector of $\bm{B}'$ is $\bm{P}^{-1}\bm{v}$ given $\bm{B}=\lambda\bm{v}$, where $\bm{P}$ is the permutation matrix from the index map given by $\psi_r$ to the index map given by $\psi'_r$.
\end{theorem}
\textbf{Proof:}
From Eq.~\eqref{eq3:sec: Row/Column Same Mapping}, we have 
\begin{eqnarray}\label{eq1: thm: unfolding maps spectrum and eigenvectors}
\bm{B}' =\bm{P}^{-1}\bm{B}\bm{P} \Longleftrightarrow \bm{P}\bm{B}'\bm{P}^{-1}=\bm{B}.
\end{eqnarray}
If $\bm{B} \bm{v} = \lambda \bm{v}$, then we have
\begin{eqnarray}\label{eq2: thm: unfolding maps spectrum and eigenvectors}
\bm{P}\bm{B}'\bm{P}^{-1}\bm{v}=\lambda \bm{v} \Longleftrightarrow \bm{B}'\bm{P}^{-1}\bm{v}=\lambda\bm{P}^{-1}\bm{v}.
\end{eqnarray}
If $\bm{v}$ is an eigenvector of $\bm{B}$ corresponding to the eigenvalue $\lambda$, then $\bm{P}^{-1}\bm{v}$ is an eigenvector of $\bm{B}'$ with the same eigenvalue. By interchanging the roles of $\bm{B}$ and $\bm{B}'$, $\bm{B}'$ and $\bm{B}$ have the same eigenvalues and the eignevector of $\bm{B}'$ is $\bm{P}^{-1}\bm{v}$ given $\bm{B}=\lambda\bm{v}$.
$\hfill \Box$

Let us present an example to illustrate Theorem~\ref{thm: unfolding maps spectrum and eigenvectors}.
\begin{example}\label{exp1: Row/Column Same Mapping}
Given a tensor $\mathcal{B}  \in \mathfrak{R}^{2 \times 2 \times 2 \times 2}$ expressed by 
\begin{eqnarray}\label{eq1: exp1}
\mathcal{B}&=&\left[
    \begin{array}{cc : cc}
       b_{1,1;1,1} & b_{1,2;1,1} & b_{1,1;1,2} & b_{1,2;1,2}   \\
       b_{2,1;1,1} & b_{2,2;1,1} & b_{2,1;1,2} & b_{2,2;1,2}   \\ \hdashline[2pt/2pt]
       b_{1,1;2,1} & b_{1,1;2,1} & b_{1,1;2,2} & b_{1,2;2,2}   \\ 
        b_{2,1;2,1} & b_{2,2;2,1} & b_{2,1;2,2} & b_{2,2;,2,2}   \\ 
    \end{array}
\right]=\left[
    \begin{array}{cc : cc}
       3 & 2 & 4 & 5  \\
       1 & 7 & 8 & 9  \\ \hdashline[2pt/2pt]
       1 & 2 & 5 & 5  \\ 
       3 & 4 & 0.1 & 9.9  \\
    \end{array}
\right]. 
\end{eqnarray}

If we apply the following unfolding map to $\mathcal{B}$
\begin{eqnarray}\label{eq2: exp1}
\psi_r(1,1,2,2)=1, \psi_r(1,2,2,2)=3, \psi_r(2,1,2,2)=2, \psi_r(2,2,2,2)=4, 
\end{eqnarray}
we have
\begin{eqnarray}\label{eq3: exp1}
\bm{B}&=&\left[
    \begin{array}{cccc}
       3 & 1 & 4 & 5  \\
       1 & 3 & 8 & 0.1  \\ 
       2 & 2 & 5 & 5  \\ 
       7 & 4 & 9 & 9.9  \\
    \end{array}
\right]. 
\end{eqnarray}
On the other hand, if we apply the following unfolding map to $\mathcal{B}$
\begin{eqnarray}\label{eq4: exp1}
\psi'_r(1,1,2,2)=2, \psi'_r(1,2,2,2)=3, \psi'_r(2,1,2,2)=4, \psi'_r(2,2,2,2)=1, 
\end{eqnarray}
we have
\begin{eqnarray}\label{eq5: exp1}
\bm{B}'&=&\left[
    \begin{array}{cccc}
       9.9 & 7 & 9 & 4  \\
       5 & 3 & 4 & 1  \\ 
       5 & 2 & 5 & 2  \\ 
       0.1 & 1 & 8 & 3  \\
    \end{array}
\right]. 
\end{eqnarray}

Then, we have the following permutation matrix $\bm{P}$:
\begin{eqnarray}\label{eq6: exp1}
\bm{P}&=&\left[
    \begin{array}{cccc}
       0 & 1 & 0 & 0  \\
       0 & 0 & 0 & 1  \\
       0 & 0 & 1 & 0  \\
       1 & 0 & 0 & 0  \\
    \end{array}
\right],
\end{eqnarray}
such that 
\begin{eqnarray}\label{eq7: exp1}
\bm{B}'&=&\bm{P}^{-1}\bm{B}\bm{P}.
\end{eqnarray}
The unfolded matrices $\bm{B}'$ and $\bm{B}$ have the same eigenvalues: $[18.57, 3.12,  -0.39+1.22\sqrt{-1}, 
 -0.39-1.22\sqrt{-1}]$. Moreover, the eignevector of $\bm{B}'$ is $\bm{P}^{-1}\bm{v}$ given $\bm{B}=\lambda\bm{v}$. For example, if $\bm{v}= [-0.37, -0.23,$
$-0.39, -0.81]^{T}$ is the eigenvector of the matrix $\bm{B}$ corresponding to the eigenvalue of $18.57$, we have $[-0.81, -0.37, -0.39, -0.23]^{T}$ is the eigenvector of the matrix $\bm{B}'$ corresponding to the eigenvalue of $18.57$. Note that $\bm{P}^{-1} [-0.37, -0.23, -0.39, -0.81]^{T} =[-0.81, -0.37, -0.39, -0.23]^{T}$.
\end{example}

\begin{remark}\label{rmk: count of permutation similarity classes same}
The count of permutation similarity classes given by Eq.~\eqref{eq7: exp1} for permutation matrices of a tensor $\mathcal{B} \in \mathfrak{R}^{I_1 \times \cdots \times I_M \times I_1 \times \cdots \times I_M}$ corresponds to the partition number $\mbox{Per}\left(\mathbb{I}_M\right)$. The partition number $\mbox{Per}\left(\mathbb{I}_M\right)$ can be expressed by
\begin{eqnarray}\label{eq8: exp1}
\mbox{Per}\left(\mathbb{I}_M\right) \sim \frac{1}{4\sqrt{3}\mathbb{I}_M}\exp\left(\sqrt{\frac{2\pi^2\mathbb{I}_M}{3}}\right),
\end{eqnarray}
given that $\mathbb{I}_M$ becomes large. Eq.~\eqref{eq8: exp1} is derived from~\cite{li2018note}.
\end{remark}

\section{Row/Column Different Mapping}\label{sec: Row/Column Different Mapping}

If we have $\psi_r\neq\psi_c$, any two unfolded matrices $\bm{B}$ and $\bm{B}'$ can be related by the following
\begin{eqnarray}\label{eq3-1: Row/Column Different Mapping}
\bm{B}' = \bm{Q}^{-1}\bm{B}\bm{P},
\end{eqnarray} 
where $\bm{P}$ and $\bm{Q}$ are two different permutation matrices, i.e., entries are $0$ and $1$ only. In this section, we will consider more general entries in invertible matrices $\bm{P}$ and $\bm{Q}$ by considering entries from $\mathfrak{R}$.

\subsection{The  Local-Global Principle of Smith Normal Form}\label{subsec: The  Local-Global Principle of Smith Normal Form}

Let $\mathfrak{R}$ be a B{\'e}zout ring with quotine field $\mathfrak{R} \subset \mathbb{F}$, and let $p$ be a fixed prime element of $\mathfrak{R}$. Then, every element $f \in \mathbb{F}$ can be expressed by the following form
\begin{eqnarray}\label{eq1: Row/Column Different Mapping}
f = \frac{r_1}{r_2}p^a,
\end{eqnarray}
where $r_1, r_2 \in \mathfrak{R}$ and they are relative prime to $p$, and $a \in \mathbb{Z}$. Since the integer $a$ is determined by $f$ uniquely, we use the symbol $\mbox{ord}_p f$ to represent this $a$ value. The \emph{localization} of $\mathfrak{R}$ with respect to the prime ideal $p$, denoted by $\mathfrak{R}_p$, can be defined as 
\begin{eqnarray}\label{eq2: Row/Column Different Mapping}
\mathfrak{R}_{(p)} \define \{f \in \mathbb{F} |\mbox{ord}_p f \geq 0\}.
\end{eqnarray}
Note that $\mathfrak{R}_{(p)}$ is a principal ideal domain, and the proper ideals of $\mathfrak{R}_{(p)}$ satisfy $(p) \supset (p^2) \supset (p^3) \supset \cdots$.

Given unfolded matrices $\bm{B}$ and $\bm{B}'$, we say that $\bm{B}$ and $\bm{B}'$ is equivalent over $\mathfrak{R}$, denoted by $\bm{B} \sim \bm{B}'$, if
\begin{eqnarray}\label{eq3: Row/Column Different Mapping}
\bm{B}' = \bm{Q}^{-1}\bm{B}\bm{P},
\end{eqnarray} 
where matrix entries of $\bm{B}' , \bm{B} \in \mathfrak{R}$. On the other hand, we say that $\bm{B}$ and $\bm{B}'$ is equivalent over $\mathfrak{R}_{(p)}$, denoted by $\bm{B} \sim_p \bm{B}'$, if
\begin{eqnarray}\label{eq4: Row/Column Different Mapping}
\bm{B}' = \bm{Q}^{-1}\bm{B}\bm{P},
\end{eqnarray} 
where matrix entries of $\bm{B}' , \bm{B} \in \mathfrak{R}_{(p)}$. If $\bm{B} \in \mathfrak{R}_{(p)}^{\mathbb{I}_M \times \mathbb{I}_M}$, then $\bm{B}$ has a \emph{Smith normal form} over $\mathfrak{R}_{(p)}$ expressed by 
\begin{eqnarray}\label{eq5: Row/Column Different Mapping}
\bm{B}\sim_p \mbox{diag}\left(p^{a_1},\cdots,p^{a_r}, 0,\cdots,0\right) \define S_p(\bm{B}),
\end{eqnarray} 
where integers $0 \leq a_1 \leq a_2 \leq \cdots \leq a_r$ and $r=\mbox{rank}\bm{B}$. We have the following theorem about $S_p(\bm{B})$. 

\begin{theorem}\label{thm: Local-Global Principle}
Let $\bm{B} , \bm{B}' \in \mathfrak{R}^{\mathbb{I}_M \times \mathbb{I}_M}$ with $\mbox{rank}(\bm{B})=r$, we have $\bm{B} \sim \bm{B}'$ if and only if $\bm{B} \sim_p \bm{B}'$ for all prime ideals  $p$. Besides, we also have
\begin{eqnarray}\label{eq1: thm:Local-Global Principle}
S(\bm{B})&=&\prod\limits_{p \in P}S_p(\bm{B}), 
\end{eqnarray}
where $S(\bm{B})$ is the \emph{Smith normal form} of $\bm{B}$, and $P$ is a set of nonassociated primes.
\end{theorem}
\textbf{Proof:}
Let $p_1,\cdots,p_k$ represent the unique prime ideals whose powers appear as elementary divisors of the matrix $\bm{B}$. Then, we have 
\begin{eqnarray}\label{eq2: thm:Local-Global Principle}
S(\bm{B})&=&\mbox{diag}\left(\prod\limits_{k=1}^K p_k^{a_{1,k}},\cdots,\prod\limits_{k=1}^K p_k^{a_{r,k}}, 0,\cdots,0\right),
\end{eqnarray}
where $0  \leq a_{i,k} \leq a_{i+1,k}$ for $1 \leq i \leq r-1$ and $1 \leq k \leq K$.  

Since $p_k$ is $p_{k'}$-adic unit if $k \neq k'$, we have 
\begin{eqnarray}\label{eq3: thm:Local-Global Principle}
S_{p_k}(\bm{B})&=&\mbox{diag}\left(p_k^{a_{1,k}},\cdots,p_k^{a_{r,k}}, 0,\cdots,0\right),
\end{eqnarray}
where $1 \leq k \leq K$. Moreover, if $p \notin \{p_1,\cdots,p_K\}$, we have
\begin{eqnarray}\label{eq4: thm:Local-Global Principle}
S_{p}(\bm{B})&=&\mbox{diag}\left(\overbrace{1,\cdots,1}^{r}, 0,\cdots,0\right).
\end{eqnarray}
This Theorem is proved by product relation between Eq.~\eqref{eq2: thm:Local-Global Principle}, and Eq.~\eqref{eq3: thm:Local-Global Principle} or Eq.~\eqref{eq4: thm:Local-Global Principle}.
$\hfill \Box$

\subsection{The Number of Equivalence Classes}\label{subsec: The Number of Equivalence Classes}

Let $\sigma(\mathfrak{R})$ be the set of $\sim$-equivalanece classes for the matrix $\bm{B} \in \mathfrak{R}^{\mathbb{I}_M \times \mathbb{I}_M}$ provided Eq.~\eqref{eq3: Row/Column Different Mapping}. We use $\tilde{\bm{B}}$ to represent the equivalence class of the matrix $\bm{B}$ and define the following addition operation with respect to any two equivalence classes $\tilde{\bm{B}}$ and $\tilde{\bm{D}}$ as
\begin{eqnarray}\label{eq:def eq classes addi}
\tilde{\bm{B}}+\tilde{\bm{D}} \define \widetilde{B \oplus D},
\end{eqnarray}
where $\oplus$ is the direct sum operation. Similarly, we also define the following product operation with respect to any two equivalence classes $\tilde{\bm{B}}$ and $\tilde{\bm{D}}$ as
\begin{eqnarray}\label{eq:def eq classes prod}
\tilde{\bm{B}}\tilde{\bm{D}} \define \widetilde{B \otimes D},
\end{eqnarray}
where $\otimes$ is the Kronecker product of two matrices. We have the following theorem to characterize $(\mathfrak{R})$ via polynomials.

\begin{theorem}\label{thm: inj semiring homo}
Let $\sigma(\mathfrak{R})$ be the set of $\sim$-equivalanece classes for the matrix $\bm{B} \in \mathfrak{R}^{\mathbb{I}_M \times \mathbb{I}_M}$. There exists an injective semiring homomorphisim, denoted by $\Psi$, which is given by
\begin{eqnarray}\label{eq1: thm: inj semiring homo}
\Psi: \sigma(\mathfrak{R}) \rightarrow \prod\limits_{p \in P}\mathbb{Z}^{+}[x_p],
\end{eqnarray}
where $\mathbb{Z}^{+}[x_p]$ is the semiring of polynomials in $x_p$ with nonnegative integral coefficients.
\end{theorem}
\textbf{Proof:}
We will apply the localization technique to $\sigma(\mathfrak{R})$. Considering $\bm{B}$ with the following Smith normal form 
\begin{eqnarray}\label{eq2: thm: inj semiring homo}
S_p(\bm{B})&=&\mbox{diag}\left(\overbrace{1,\cdots,1}^{n_{p,0}},\overbrace{p,\cdots,p}^{n_{p,1}},\cdots
\overbrace{p^{r_p},\cdots,p^{r_p}}^{n_{p,r_p}}\right),
\end{eqnarray}
then, we define 
\begin{eqnarray}\label{eq3: thm: inj semiring homo}
\Psi_p (\tilde{\bm{B}})\define\sum\limits_{k=0}^{r_p}n_{p,k}x_p^k,
\end{eqnarray}
where $n_{p,k}$ are nonnegative integers. Note that R.H.S. of Eq.~\eqref{eq3: thm: inj semiring homo} belongs to $\mathbb{Z}^{+}[x_p]$. From Eq.~\eqref{eq:def eq classes addi} and Eq.~\eqref{eq:def eq classes prod}, the map provided by Eq.~\eqref{eq3: thm: inj semiring homo} is a homomorphism map since
\begin{eqnarray}\label{eq4: thm: inj semiring homo}
\bm{B}\oplus\bm{D} &\sim_p& S_p(\bm{B}) \oplus S_p(\bm{D}),\nonumber \\
\bm{B}\otimes\bm{D} &\sim_p& S_p(\bm{B}) \otimes S_p(\bm{D}).
\end{eqnarray}

By setting 
\begin{eqnarray}\label{eq5: thm: inj semiring homo}
\Psi(\tilde{\bm{B}}) \define \prod\limits_{p \in P}\Psi_p (\tilde{\bm{B}}).
\end{eqnarray}
We have $\Psi$ is a homomorphism, and injective by Theorem~\ref{thm: Local-Global Principle}. This theorem is proved.  
$\hfill \Box$

From Theorem~\ref{thm: inj semiring homo}, we can have the following corollary to determine the number of equivalence classes in the set $\sigma(\mathfrak{R})$, represented by $|\sigma(\mathfrak{R})|$.
\begin{corollary}\label{cor: class number bounds}
Let $\sigma(\mathfrak{R})$ be the set of $\sim$-equivalanece classes for the matrix $\bm{B} \in \mathfrak{R}^{\mathbb{I}_M \times \mathbb{I}_M}$. Then, we have
\begin{eqnarray}\label{eq1: cor: class number bounds}
|\sigma(\mathfrak{R})|=\left(\sum\limits_{k=0}^{\mathbb{I}_M}H_{\mathbb{I}_M - 1}^{k+1}\right)^{|P|},
\end{eqnarray}
where $|P|$ is the number of distinct nonassociated primes in $P$. 
\end{corollary}
\textbf{Proof:}
According to Eq.~\eqref{eq2: thm: inj semiring homo}, we have
\begin{eqnarray}\label{eq2: cor: class number bounds}
n_{p,0}+n_{p,1}+\cdots+n_{p,r_p} = \mathbb{I}_M.
\end{eqnarray}
For $r_p=k$, Eq.~\eqref{eq2: cor: class number bounds} can be expressed by
\begin{eqnarray}\label{eq3: cor: class number bounds}
n_{p,0}+n_{p,1}+\cdots+n'_{p,k} = \mathbb{I}_M-1,
\end{eqnarray}
where $n_{p,0} \geq 0$, $n_{p,1} \geq 0$, $\cdots$, $n'_{p,k} \geq 0$, and $0 \leq k \leq \mathbb{I}_M$. By applying 
formula of combination with repetitions to Eq.~\eqref{eq3: cor: class number bounds}, there are $H_{\mathbb{I}_M - 1}^{k+1}$ solutions of $n_{p,0},n_{p,1},\cdots,n'_{p,k}$ in Eq.~\eqref{eq3: cor: class number bounds}, i.e., there are 
 $H_{\mathbb{I}_M - 1}^{k+1}$ different $\Psi_p (\tilde{\bm{B}})$ given $r_p=k$. 

Because $0 \leq k \leq \mathbb{I}_M$, the number of different $\Psi_p (\tilde{\bm{B}})$ can be determined by
\begin{eqnarray}\label{eq4: cor: class number bounds}
\sum\limits_{k=0}^{\mathbb{I}_M}H_{\mathbb{I}_M - 1}^{k+1}. 
\end{eqnarray}
Finally, from Theorem~\ref{thm: inj semiring homo} (Eq.~\eqref{eq5: thm: inj semiring homo}), this corollary is proved. 
$\hfill \Box$

Several remarks will be provided before the end of this section.
\begin{remark}\label{rmk: count of permutation similarity classes diff}
According to Eq.~\eqref{eq8: exp1}, the number of different permutation similarity classes will be $\mbox{Per}\left(\mathbb{I}_M\right)^2$ since $\bm{P}$ may not be equal to $\bm{Q}$, where both $\bm{P}$ and $\bm{Q}$ are permutation matrices. This can be approximated by 
\begin{eqnarray}\label{eq1: Remark 2}
\mbox{Per}\left(\mathbb{I}_M\right)^2 \sim \frac{1}{48\mathbb{I}^2_M}\exp\left(2\sqrt{\frac{2\pi^2\mathbb{I}_M}{3}}\right),
\end{eqnarray}
given that $\mathbb{I}_M$ becomes large.
\end{remark}

\begin{remark}\label{rmk: count of permutation similarity classes diff}
If $\mathfrak{R}$ is a field, Corollary~\ref{cor: class number bounds} becomes $|\sigma(\mathfrak{R})|=\mathbb{I}_M$ since $\sigma(\mathfrak{R})$ is isomorphism to natural number. 
\end{remark}

\begin{remark}\label{rmk: eigenvector shift diff}
If $\bm{B}=\bm{Q}^{-1}\bm{A}\bm{P}$ (equivalent class) and $\bm{A}\bm{v}=\lambda\bm{v}$, we have $\bm{B}\bm{P}^{-1}\bm{Q}(\bm{Q}^{-1}\bm{v}) = \lambda(\bm{Q}^{-1}\bm{v})$, same eigenvalues, for different matrix $\bm{B}\bm{P}^{-1}\bm{Q}$ and eigenvectors $(\bm{Q}^{-1}\bm{v})$
\end{remark}

\begin{remark}\label{rmk: eigenvector shift diff}
Theorem~\ref{thm: Local-Global Principle} and Theorem~\ref{thm: inj semiring homo} generalize the matrix equivalence class domain over PID ring discussed in~\cite{gerstein1977local} to the matrix equivalence class domain over B{\'e}zout ring. According to Theorem 2.1 from~\cite{stanley2016smith}, a B{\'e}zout ring also has Smith normal form. Therefore, the localization technique adopted by~\cite{gerstein1977local} based on Smith normal form over ring with PID can be applied to Smith normal form over B{\'e}zout ring. 
\end{remark}

\bibliographystyle{IEEETran}
\bibliography{TenFoldingByMap_Bib}

\end{document}